\documentclass{amsart}
\usepackage{graphicx}
\usepackage{float}
\newtheorem{thm}{Theorem}[section]
\newtheorem{cor}[thm]{Corollary}
\newtheorem{lem}[thm]{Lemma}
\newtheorem{prop}[thm]{Proposition}
\theoremstyle{definition}

\newtheorem{exa}[thm]{Example}
\theoremstyle{remark}

\numberwithin{equation}{section}
\begin{document}

\title{ nilpotent Lie algebras having the Schur multiplier of maximum dimension }%
\author{Afsaneh shamsaki}%
\email{Shamsaki.Afsaneh@yahoo.com}%
\author{Peyman Niroomand}%
\email{p$\_$niroomand@yahoo.com}
\address{Damghan university}%
\email{Shamsaki.Afsaneh@yahoo.com}%

%\thanks{0}%
%\subjclass{0}%
\keywords{ nilpotent Lie algebra}%

%\date{0}%
%\dedicatory{0}%
%\commby{0}%
% ----------------------------------------------------------------
\begin{abstract}
Let $ L $ be an $ n $-dimensional nilpotent Lie algebra of nilpotency class $ c $ with the derived subalgebra of dimension $ m $. Recently,  Rai proved that the dimension of  Schur multiplier of $ L $ is bounded by $ \frac{1}{2}(n-m-1)(n+m)-\sum\limits_{i=2}^ {min\lbrace n-m,c\rbrace} n-m-i $. In this paper, we obtain the structure of all nilpotent Lie algebras that attain this bound.
 \end{abstract}
\maketitle
% ----------------------------------------------------------------
\section{Introduction }
The Schur multiplier of groups was introduced by Schur in \cite{9} on the study of projective representation of groups.  It is known from \cite[Main Theorem]{4} that for a non-abelian $ p $-group $ G $ of order $p^{n}$ and the derived  subgroup of order $ p^{k} $, we have
\begin{equation}\label{e1-}
\mid \mathcal{M}(G) \mid \leq p^{\frac{1}{2}(n-k-1)(n+k-2)+1}.
\end{equation}
Recently Rai in \cite[Theorem 1.1]{6-} classified $ p $-groups of nilpotency class two such that $ \mid \mathcal{M}(G) \mid = p^{\frac{1}{2}(n-k-1)(n+k-2)+1}$ and then Hatui in \cite{1} showed for $ p \neq 3 $ only nilpotent groups of class 2 satisfy in the upper bound \ref{e1-}. By  a classical result due to Lazard  \cite{2-}, we may associate a $ p $-group to a Lie algebra. Analogous to the Schur multiplier of a group, the Schur multiplier of a Lie algebra, $ \mathcal{M}(L) $, can be defined as $ \mathcal{M}(L)\cong R\cap F^{2}/ [R, F] $ where $ L\cong F/R $ and $ F $ is a free Lie algebra by the same motivation in \cite[Theorem 3.1]{3}, the second author  proved that for a non-abelian nilpotent Lie algebra of dimension $ n $ and the derived algebras of dimension $ m $, we showed
\begin{equation}\label{e1}
\dim \mathcal{M}(L)\leq \frac{1}{2}(n+m-2)(n-m-1)+1.
\end{equation}
In particular, when $ m=1 $ the bound is attained if and only if $ L\cong H(1)\oplus A(n-3) $ where $ H(m) $ and $ A(n) $ are used to denote the Heisenberg and Abelian Lie algebras of dimension $ 2m+1 $ and $ n $, respectively.
Looking \cite{5} shows similar results of  Hatui are obtained for the class of nilpotent Lie algebra.
Recently, Rai in \cite[Theorem 1.1]{7} improved the bound \ref{e1} and showed
\begin{equation}\label{e2}
\dim \mathcal{M}(L)\leq \frac{1}{2}(n-m-1)(n+m)-\sum\limits_{i=2}^ {min\lbrace n-m,c\rbrace} n-m-i,
\end{equation}
for an $ n $-dimensional nilpotent Lie algebra $ L $ of nilpotency class $ c $ when $ \dim L^{2}=m\geq 1 $.\\
In the current paper, we are going to classify all nilpotent Lie algebras that attain the bound \ref{e2}. More precisely, we show that such Lie algebras are nilpotent of class two.\\
\section{Preliminaries}
For the convenience of the reader, we give some results without proofs which will be used in the next section.\\
We recall two theorems.
\begin{thm}$($See \cite[Corollary 2.3]{3}$)$\label{2}
Let $L$ be a finite dimensional Lie algebra and $K$ be a central ideal of $L$. Then
\begin{align*}
\dim\mathcal{M}(L)+\dim(L^{2}\cap K) \leq \dim\mathcal{M}(L/K)+\dim\mathcal{M}(K)+\dim ((L/K)^{ab}\otimes_{ab} K) .
\end{align*}
\end{thm}
\begin{thm}$($See \cite[Theorem 3.1]{4}$)$\label{3}
Let $ L $ be an $ n $-dimensional nilpotent Lie algebra $ L $, $ \dim L^{2}=n-2 $ and $ n\geq 4 $. Then $ \dim \mathcal{M}(L)\leq \dim L^{2} $.
\end{thm}
Let $ Z^{\wedge}(L) $ is used to denote the exterior center of a Lie algebra $ L $ (See \cite{2--} for more information). By \cite{2--}, $ L $ is capable if and only if $ Z^{\wedge}(L)=0 $.\\
Here is easy result.
\begin{lem}$($See \cite[Corollary 2.3]{6}$)$\label{3-}
Let $ L $ be a finite dimensional non-abelian nilpotent Lie algebra. Then $ Z^{\wedge}(L)\subseteq L^{2} $.
\end{lem}
Let $ \otimes_{ab} $ denotes the operator of the standard tensor product of Lie algebras (See \cite[page $103$]{1-}). Then
\begin{prop}$($See \cite[Proposition 2.3]{5}$)$\label{4}
Let $ L $ be a Lie algebra. Then
\begin{itemize}
\item[(i)] the map $ \gamma _{L} : L^{ab} \otimes _{mod} L^{ab} \otimes _{mod} L^{ab} \rightarrow L^{2}/ L^{3} \otimes _{mod} L/L^{2} $ given by $$ (x+ L^{2}) \otimes (y+ L^{2}) \otimes (z+ L^{2}) \mapsto ([x, y]+L^{3}\otimes z+L^{2})+([z,x]+L^{3}\otimes y+L^{2})+([y,z]+L^{3}\otimes x+L^{2})$$ is a Lie homomorphism. If any two element of the set $ \lbrace x, y, z\rbrace \subseteq L $ are linearly dependent. Then $ \gamma _{L}(x+ L^{2} \otimes y+ L^{2} \otimes z+ L^{2})=0$.
\item[(ii)] Define the map
$$ \gamma ^{\prime} _{2}:( L/ Z(L))^{ab} \otimes _{mod} (L/Z(L)) ^{ab} \otimes _{mod} (L/Z(L)) ^{ab} \rightarrow L^{2}/ L^{3} \otimes _{mod} (L/Z(L))^{ab} $$
\begin{align*}
&(x+ (L^{2}+Z(L))) \otimes (y+ (L^{2}+Z(L))) \otimes (z+ (L^{2}+Z(L))) \otimes (w+ (L^{2}+Z(L)))\mapsto \cr
&([x, y]+L^{3}\otimes z+(L^{2}+Z(L)))+([z,x]+L^{3}\otimes y+(L^{2}+Z(L)))+ \cr
&([y,z]+L^{3}\otimes x+(L^{2}+Z(L))).
\end{align*}
Then $\gamma ^{\prime} _{2}$ is a Lie homomorphism. Moreover, if any two element of the set $ \lbrace x, y, z\rbrace \subseteq L$ are linearly dependent. Then $ \gamma ^{\prime} _{2}(x+ L^{2} \otimes y+ L^{2} \otimes z+ L^{2})=0$.
\item[(iii)] The map
$\gamma ^{\prime} _{3}:( L/ Z(L))^{ab} \otimes _{mod} (L/Z(L)) ^{ab} \otimes _{mod} (L/Z(L)) ^{ab} \otimes _{mod} (L/Z(L)) ^{ab} $ \\
$\rightarrow L^{3} \otimes _{mod} (L/Z(L))^{ab} $ given by
\begin{align*}
&(x+ (L^{2}+Z(L))) \otimes (y+ (L^{2}+Z(L))) \otimes (z+ (L^{2}+Z(L))) \otimes (w+ (L^{2}+Z(L))) \cr
& \mapsto ([[x, y],z] \otimes w+(L^{2}+Z(L)))+([w,[x,y]] \otimes z+(L^{2}+Z(L)))+ \cr
&([[z,w],x]\otimes y+(L^{2}+Z(L)))+([y,[z,w]] \otimes x+(L^{2}+Z(L))).
\end{align*}
\end{itemize}
is a Lie homomorphism.
\end{prop}
Now there are some results on the dimension of the exterior square of a Lie algebra $ L $.
\begin{thm}$($See \cite[Theorem 2.7]{5}$)$ \label{T41}
Let $ L $ be a finite dimensional nilpotent non-abelian Lie algebra of class $ c $. Then
\begin{align*}
\dim L\wedge L + \dim \mathrm{Im} \gamma ^{\prime} _{2} &\leq \dim L\wedge L+ \sum\limits_{i=2}^ {c} \dim \ker \alpha _{i}\cr
& =\dim L/ L^{2} \wedge L/ L^{2}+ \sum\limits_{i=2}^ {c} \dim (L^{i}/L^{i+1}\otimes _{mod} (L/Z(L))^{ab}).
\end{align*}
\end{thm}
\begin{thm}$($See \cite[Theorem 2.8]{5}$)$ \label{T42}
Let $ L $ be a finite dimensional nilpotent non-abelian Lie algebra of class $ 3 $. Then, we have
\begin{align*}
&\dim L\wedge L + \dim \mathrm{Im} \gamma ^{\prime} _{2}+ \dim \mathrm{Im} \gamma ^{\prime} _{3} \cr
& \leq \dim L/ L^{2} \wedge L/ L^{2}+ \dim (L^{2}/L^{3}\otimes _{mod} L^{ab})+ \dim (L^{3}\otimes _{mod} L^{ab}).
\end{align*}
\end{thm}
Let $ Z^{*}(L) $ be a  symbol, the epicenter of a Lie algebra $ L $. From \cite[Lemma 3.1]{2--}, $ Z^{*}(L)=Z^{\wedge}(L) $.
\begin{thm}$($See \cite[Theorem 4.4]{10}$)$\label{4-}
Let $ N $ be a central ideal of a Lie algebra $ L $. Then the following condition are equivalent.
\item[(i)] $ \mathcal{M}(L)\cong \mathcal{M}(L/N)/N\cap L^{2}$ ;
\item[(ii)] $ N\subseteq Z^{*}(L) =Z^{\wedge}(L)$;
\item[(iii)] The natural map $ \mathcal{M}(L) \longrightarrow \mathcal{M}(L/N)$ is monomorphism.
\end{thm}
Let $ cl(L) $ is used to denote the nilpotency class of $ L $.\\
The following technical results are suitable for the next investigations.
\begin{thm}$($See \cite[Theorem 5.5]{6}$)$ \label{5}
Let $L$ be an $n$-dimensional nilpotent Lie algebra such that $\dim L^{2} \leq 2$. Then $ L $ is capable if and only if $ L $ is isomorphic to one the following Lie algebras.
\begin{itemize}
\item[(i)] If $ \dim L^{2}=0 $, then $ L\cong A(n) $ and $ n\geq 2 $;
\item[(ii)]If $ \dim L^{2}=1 $, then $ L\cong H(1)\oplus A(n-3) $;
\item[(iii)]If $ \dim L^{2}=2 $ and $ cl(L)=2 $, then $ L\cong L_{6,7}^{2}(\eta)\oplus A(n-6) $, $ L\cong L_{5.8} \oplus A(n-5) $, $ L\cong L_{6,22}(\varepsilon)\oplus A(n-6) $, or $ L\cong L_{1} \oplus A(n-7) $;
\item[(iv)] If $ \dim L^{2}=2 $ and $ cl(L)=3 $, then $ L\cong L_{4,3}\oplus A(n-4) $, or $ L\cong L_{5,5}\oplus A(n-5) $.
\end{itemize}
\end{thm}
\begin{thm}$($See \cite[Theorem 3.2]{8}$)$ \label{6}
There is no Lie algebra $ L $ of dimension $ n $ such $ \dim L^{2}=3 $ and $ \dim \mathcal{M}(L)=1/2(n-1)(n-2)-2 $.
\end{thm}
\begin{thm}$($See \cite[Theorem 2.22]{5}$)$ \label{7}
Let $ L $ be an $ n $-dimensional nilpotent Lie algebra of class two with the drived subalgebra of dimension $ m $. Then $ \dim \mathcal{M}(L)= \frac{1}{2}(n+m-2)(n-m-1)+1 $ if and only if $ L $ is isomorphic to one of the Lie algebras $ H(1)\oplus A(n-3) $, $ L_{5,8} $ or $ L_{6,26} $.
\end{thm}
\section{main result}
In this section, we classify all nilpotent Lie algebras $ L $ such that obtain the upper bound \ref{e2}. Throughout the paper, we say that $ \dim \mathcal{M}(L) $ attains the upper  bound provided that $ \dim \mathcal{M}(L) $ attains the upper bound \ref{e2}.
\begin{lem}\label{21}
There is no $ n $-dimensional nilpotent Lie algebra with the derived subalgebra of dimension $ n-2 $ $( n\geq 4 )$ such that $ \dim \mathcal{M}(L) $ attains the bound.
\begin{proof}
Let $ \dim \mathcal{M}(L) $ attains the bound then $ c\leq n-1 $ and  $ \dim \mathcal{M}(L)=\frac{1}{2}(n-(n-2)-1)(n+(n-2))-\sum\limits_{i=2}^ {min\lbrace 2,c\rbrace} n-m-i =n-1 $. On the other hand, Theorem \ref{3} shows $ \dim \mathcal{M}(L) \leq n-2$.  It is a contradiction.
\end{proof}
\end{lem}
\begin{thm}\label{22}
Let $ L $ be an $ n $-dimensional nilpotent Lie algebra of nilpotency class two and $ \dim L^{2}=m $. Then   $\dim \mathcal{M}(L) $ attains the bound if and only if $ L $ is isomorphic to one of the Lie algebras $ H(1)\oplus A(n-3) $, $ L_{5,8} $ or $ L_{6,26} $.
\begin{proof}
 It is obvious that the bounds \ref{e1} and \ref{e2} are equal when $ c=2 $. Since $\dim \mathcal{M}(L) $  attains the bound  and $ L $ is of nilpotency class two,  $ L $ is isomorphic to  one of the Lie algebras $ H(1)\oplus A(n-3) $, $ L_{5,8} $ or $ L_{6,26} $ by using Theorem \ref{7}.
\end{proof}
\end{thm}
\begin{prop}\label{23}
Let $ L $ be an $ n $-dimensional nilpotent Lie algebra of nilpotency class $ c\geq 3 $  such that   $\dim \mathcal{M}(L) $ attains the bound. If $ I $ is an  ideal of one-dimensional  contained in $ Z(L) \cap L^{2} $.   Then $ \dim\mathcal{M}(L/I) $ attains the bound.
\begin{proof}
Let $ I $ be an ideal of one-dimensional  such that contained in $ Z(L) \cap L^{2} $, then $ \dim L/I=n-1 $, $ \dim (L/I)^{2}=m-1 $ and $ cl(L/I)=c $ or $ c-1 $.\\
Let $ cl(L/I)=c $. Then by using  Theorem   \ref{2}, we have
\begin{align*}
\dim\mathcal{M}(L)+ \dim(L^{2}\cap I) & \leq \dim\mathcal{M}(L/I)+\dim\mathcal{M}(I)+\dim ((L/I)^{ab}\otimes I) \cr
&=\dim\mathcal{M}(L/I)+\dim ((L)^{ab}\otimes I).
\end{align*}
Since $ cl(L/I)=c $, (\ref{e2}) implies that
\begin{align*}
\dim\mathcal{M}(L) & \leq \dim\mathcal{M}(L/I)+\dim ((L)^{ab}\otimes I)- \dim(L^{2}\cap I) \cr
&\leq \frac{1}{2}(n-m-1) (n+m-2)-(\sum\limits_{i=2}^ {min\lbrace n-m,c\rbrace} n-m-i) +(n-m)-1 \cr
&= \frac{1}{2}(n-m-1)(n+m)-\sum\limits_{i=2}^ {min\lbrace n-m,c\rbrace} n-m-i= \dim \mathcal{M}(L).
\end{align*}
Therefore
\begin{equation*}
\dim \mathcal{M}(L/I)=\frac{1}{2}(n-m-1)(n+m-2)-\sum\limits_{i=2}^ {min\lbrace n-m,c \rbrace} n-m-i.
\end{equation*}
If $cl(L/I)= c-1 $. Then let $ n-m\leq c $. So we can see that $ \sum\limits_{i=2}^ {min\lbrace n-m,c \rbrace} n-m-i=\sum\limits_{i=2}^ {min\lbrace n-m,c-1 \rbrace} n-m-i $. Hence by  a similar  way to the above case,  the result follows.\\
Suppose that $ cl(L/I)=c-1 $ and $ n-m> c $. Since $ cl(L/I) = c-1 $, we have $ I= \gamma _{c}(L)$ and $ \dim \gamma _{c}(L)=1 $. On the other hand, by using the proof of \cite[Theorem 1.1]{7}, we have $ \dim \ker(\lambda_{c}) \geq n-m-c $ and
\begin{equation*}
\dim \mathcal{M}(L)=\dim \mathcal{M}(L/\gamma_{c}(L))+ \dim(L/\gamma_{2}(L)-1) \dim \gamma_{c}(L) - \dim \ker(\lambda_{c}).
\end{equation*}
Thus
\begin{align*}
& \dim\mathcal{M}(L) = \dim \mathcal{M}(L/I)+( \dim L/\gamma_{2}(L)-1) \dim I - \dim\ker(\lambda_{c}) \cr
& \leq \frac{1}{2}(n-m-1) (n+m-2)-(\sum\limits_{i=2}^ {c-1} n-m-i) +(n-m-1)-(n-m-c) \cr
&= \dim \mathcal{M}(L),
\end{align*}
and so $ \dim \mathcal{M}(L/I)= \frac{1}{2}(n-m-1)(n+m-2)-(\sum\limits_{i=2}^ {min\lbrace n-m,c\rbrace} n-m-i)$.\\
This completes the proof.
\end{proof}
\end{prop}

\begin{lem}\label{24}
Let $ L $ be an $ n$-dimensional nilpotent Lie algebra of nilpotency class $ c\geq 3 $ and $ \dim L^{2}=m\geq 2 $. If $ \dim \mathcal{M}(L) $ attains the bound, then $ L $ is capable.
\begin{proof}
By contrary, let $ L $ be  non-capable. Then by using Lemma \ref{4-} $ Z^{\wedge}(L)\neq 0 $ and Lemma \ref{3-} shows that  there is an ideal $ I $ in $ Z^{\wedge}(L) $  of dimension one such that $ I\subseteq  Z(L) \cap L^{2}$.  Theorem \ref{4-} shows that $ \dim \mathcal{M}(L)= \dim \mathcal{M}(L/I)-1$.  Using  Proposition \ref{23}, we have $ \dim \mathcal{M}(L/I)=\frac{1}{2}(n-m-1)(n+m-2)-(\sum\limits_{i=2}^ {min\lbrace n-m,c^{\prime}\rbrace} n-m-i) $, in which $ c^{\prime}=cl(L/I) $. Since $ c^{\prime}=c $ or $ c-1 $, we have $ \dim \mathcal{M}(L/I)= \frac{1}{2}(n-m-1)(n+m-2)-(\sum\limits_{i=2}^ {min\lbrace n-m,c^{\prime}\rbrace} n-m-i) $, and so $\dim \mathcal{M}(L)=  \dim \mathcal{M}(L/I)-1= \frac{1}{2}(n-m-1)(n+m-2)-(\sum\limits_{i=2}^ {min\lbrace n-m,c\rbrace} n-m-i)-1$ or $ \frac{1}{2}(n-m-1)(n+m-2)-(\sum\limits_{i=2}^ {min\lbrace n-m,c-1\rbrace} n-m-i) $. It is a contradiction.
\end{proof}
\end{lem}
\begin{lem}\label{25}
There is no $ n $-dimensional nilpotent Lie algebra $ L $ of nilpotency class $ 3 $  such that $ \dim  \mathcal{M}(L) $ attains the bound and $ cl(L/I) = 2 $ for a given one-dimensional ideal $ I $ in $  Z(L)\cap L^{2} $.
\begin{proof}
By contrary, let $ L $ be an $ n $-dimensional nilpotent Lie algebra of nilpotency class
3  such that $ \dim \mathcal{M}(L) $ attains the bound. Let $ cl(L/I) = 2 $ for a given one-dimensional ideal $ I $ in $  Z(L)\cap L^{2} $  . By Proposition \ref{23}, $ \dim \mathcal{M}(L/I)$ also attains the bound. Hence by Theorem \ref{22},  $ L/I $  is isomorphic to one of the Lie algebras $ H(1)\oplus A(n-3) $, $ L_{5,8} $ or $ L_{6,26} $. Now we consider all cases. Let $ \dim L^{2}=m $. Clearly $ m\geq 2 $. 

Case $(i) $. Here $ L/I\cong   H(1)\oplus A(n-4) $, then $ m=2 $. Since  $ cl(L)=3 $    and $ \dim \mathcal{M}(L) $  attains the bound. $ L $ is capable by using Lemma \ref{24}. Now $ L $ is isomorphic to one of the Lie algebras  $ L \cong L_{4,3}\oplus A(n-4)$ or $ L_{5,5}\oplus A(n-5) $ by using Theorem \ref{5}.\\
 Let $ L \cong L_{4,3}\oplus A(n-4)$, we have $ \dim \mathcal{M}(L)=2+\frac{1}{2}(n-4)(n-5)+2(n-4) $. Since $ \dim \mathcal{M}(L) $ attains the bound, we have
\begin{equation*}
\frac{1}{2}(n-3)(n+2)-(2n-9)=2+\frac{1}{2}(n-4)(n-5)+2(n-4),
\end{equation*}
which is a contradiction.\\
Let now $ L \cong L_{5,5}\oplus A(n-5)$ by a similar way to the  above case, we get  a contradiction again.\\
Case $(ii) $. Let $ L/I\cong  L_{5,8} $, then  $ n=6 $, $ m=3 $ and $ \dim \mathcal{M}(L)=8=1/2(6-1)(6-2)-2 $. It is a contradiction by looking Theorem \ref{6}.\\
Case $(iii) $. Let $ L/I\cong  L_{6,26} $ and $ L $ be minimally generated by $ d $ elements. Then   $ d=3$ and $ \dim \mathcal{M}(L)=10 $. If $ Z(L)\subseteq L^{2} $, then we can choose a basis $ B=\lbrace x_{1}+L^{2},  x_{2}+L^{2}, x_{3}+L^{2} \rbrace$ for $ L^{ab} $. We may assume that  $ [x_{1}, x_{2}]+L^{3} $ is non-trivial in $ L^{2}/L^{3} $. Since $  [x_{1}, x_{2}]\notin L^{3} $ and $ x_{2}\notin L^{2} $,  $ \gamma ^{\prime} _{2} (x_{1}+L^{2} \otimes x_{2}+L^{2} \otimes x_{3}+L^{2})\neq 0 $. Thus $ \dim \mathrm{Im} \gamma ^{\prime} _{2} \geq 1$. Let $ y \in L \setminus L ^{2} $, then $  \gamma ^{\prime} _{3} (x_{1}+L^{2} \otimes x_{2}+L^{2} \otimes x_{3}+L^{2}+ \otimes y+L^{2}) \neq 0$, hence $ \dim \mathrm{Im} \gamma ^{\prime} _{3} \geq 1$. By using Theorem \ref{T42}, we have $ \dim \mathcal{M}(L)\leq 9 $, which is a contradiction.\\
If $ Z(L)\not\subset L^{2}  $ then $ t=\dim Z(L)/ (Z(L)\cap L^{2})\geq 1 $. By using Theorem \ref{T41}, we have $ \dim \mathcal{M}(L)\leq 7 $,  which is a contradiction.\\
Hence the assumption is false and the result follows.
\end{proof}
\end{lem}
\begin{thm}\label{26}
There is no $ n $-dimensional nilpotent Lie algebra of nilpotency class $ 3 $  such that $\dim \mathcal{M}(L) $ attains the bound.
\begin{proof}
By contrary, let there be a Lie algebra $ L $ such that $ \dim \mathcal{M}(L) $ attains the bound. Consider an  ideal  $ I $ in $  Z(L)\cap L^{2} $ of dimension one. By Proposition \ref{23} $ \dim \mathcal{M}(L/I) $ also attains the bound. If $ cl(L/I) = 2 $, then we have a contradiction by using  Lemma \ref{25}. Now  if $ cl(L/I) =3$, then we proceed by induction on $ n $ to show that there is no such Lie algebra $ L $  such that $ \dim \mathcal{M}(L) $  attains the bound. For $ n=5 $, $ L/I $ is of maximal class and so $ \dim (L/I)^{2}=2 $. Hence   by using Lemma \ref{21} and Proposition \ref{23} there is no such a Lie algebras. So we have a contradiction. Let $ n> 5 $. If $ K/I $ is a one-dimensional ideal of $ L/I $ is contained in $ Z(L/I)\cap (L/I)^{2} $  such that $ (L/I)/(K/I) $  is of class $ 2 $, then the result follows from  Proposition \ref{23} and Lemma \ref{25}. Now consider $ cl((L/I)/(K/I)) =3$. Since $ \dim L/I=n-1 $, by using the induction hypothesis  there is no such a Lie algebra $ L/I $ such that  $ \dim \mathcal{M}(L/I) $ attains the bound, which is contradiction.
\end{proof}
\end{thm}

\begin{thm}\label{27}
There is no $ n $-dimensional nilpotent Lie algebra of nilpotency class $c\geq 3 $  such that $\dim \mathcal{M}(L) $ attains the  bound. That means,  $\dim \mathcal{M}(L)\leq\frac{1}{2}(n-m-1)(n+m)-(\sum\limits_{i=2}^ {min\lbrace n-m,c\rbrace} n-m-i)-1  $, for all Lie algebras of nilpotency class $ c\geq 3 $.
\begin{proof}
We proceed by using  induction on $ c $. If $ c=3 $ and $ I $ be a one dimensional ideal in $ L^{2}\cap Z(L) $. Then the result follows from Theorem \ref{26}. Now let $ c> 3 $. If  $ cl(L/I)< c $, then by using the induction hypothesis on $ c $, there is no $ L/I $ such that   $\dim \mathcal{M}(L/I)$ attains the bound. Hence  the result follows by Proposition \ref{23}. Now let $ cl(L/I)=c $.   If $K/I$ is a one dimensional  ideal  $ L/I$ such that  $ K/I\subseteq Z(L/I)\cap (L/I)^{2} $. Then, we use induction on $ n $ to prove our result. Obviously $ n\geq c+1 $ and Lemma \ref{21} emphasis $ n>c+1 $. If $ n=c+2 $, then $ L/I $ is of maximal class. So the result is obtained by using $ n=c+2 $ by using Lemma \ref{21} and Proposition \ref{23}.  Let  $ n> c+2 $.  If  $ K/I$  is an ideal in $ L/I $ of one dimensional such that $ K/I\subseteq Z(L/I)\cap (L/I)^{2} $  and $cl(L/I)/(K/I))< c $  then the result follows by using the induction hypothesis on $ c $ and by Proposition \ref{23}. Now if  $cl((L/I)/(K/I)= c $. Since $ \dim L/I=n-1 $, by the induction hypothesis on $ n $ there is no such  Lie algebra $ L/I $ that  $ \dim \mathcal{M}(L/I) $ attains the bound. Hence by using  Proposition \ref{23} is completes the proof.
\end{proof}
\end{thm}
\begin{cor}
Let $ L $ be a non-abelian nilpotent Lie algebra of class c. Then $\dim \mathcal{M}(L) $ attains the bound  if and only if $ L $ is isomorphic to one of the Lie algebras  $ H(1)\oplus A(n-3) $, $ L_{5,8} $ or $ L_{6,26} $.
\begin{proof}
The result follows from Theorems \ref{22} and \ref{27}.
\end{proof}
\end{cor}
In the following examples, we may use method of Hardy and Stitzinger in \cite{2} to show that there are some Lie algebras of dimension $ n $ with the derived subalgebra of dimension $ m $ such that $ \dim \mathcal{M}(L)= \frac{1}{2}(n-m-1)(n+m)-(\sum\limits_{i=2}^ {min\lbrace n-m,c\rbrace} n-m-i)-1 $.
\begin{exa}
Let $L\cong L_{5,7}=\langle x_{1}, ..., x_{5}\mid [x_{1}, x_{2}]=x_{3},[x_{1}, x_{3}]=x_{4}, [x_{1}, x_{4}]=x_{5} \rangle $  or $ L\cong L_{5,9}=\langle x_{1}, ..., x_{5}\mid [x_{1}, x_{2}]=x_{3},[x_{1}, x_{3}]=x_{4}, [x_{2}, x_{3}]=x_{5} \rangle $.\\
By using the method of Hardy and Stitzinger in \cite{2} we have $ \dim \mathcal{M}(L_{5,7})=\dim \mathcal{M}(L_{5,9})=3  $. Therefore $ L_{5,7} $ and $ L_{5,9} $ obtain the upper bound mentioned in the Theorem \ref{27}.
\end{exa}

\end{document}